\documentclass[a4paper,12pt]{article}
\usepackage{graphicx}
\usepackage{epsfig}

\usepackage{amsmath,amsfonts,amssymb,amsthm}
\theoremstyle{plain}
\newtheorem{thm}{Theorem}
\newtheorem*{thm*}{Theorem}

\newtheorem{lem}{Lemma}
\newtheorem*{lem*}{Lemma}

\theoremstyle{definition}

\newcommand{\reftit}{\textit}    
\newcommand{\refis}{\textbf}     

\begin{document}

\title{Multi-particle processes with reinforcements.}
\author{Yevgeniy Kovchegov 
 \footnote{This research was supported in part by NSF VIGRE Grant DMS
 9983726 at UCLA}\\
 Department of Mathematics,\\ 
 Oregon State University\\
 Corvallis, OR\\
 97331-4605, USA\\
 \texttt{kovchegy@math.oregonstate.edu}}
\date{ }
\maketitle

\begin{abstract}
 We consider a multi-particle generalization of linear edge-reinforced random walk (ERRW).
 We observe that in absence of exchangeability, new techniques are needed in order to study the multi-particle model. 
 We describe an unusual coupling construction associated with the two-point edge-reinforced process on $\mathbb{Z}$ and prove a form of recurrence: the two particles
 meet infinitely often a.s.
\end{abstract}

\section{Introduction}

  The edge-reinforced random walk was first introduced in \cite{cd} and \cite{d}. 
  In this paper we will study linear multi-particle edge-reinforced processes on $\mathbb{Z}$.  In the original
  edge-reinforced random walk model, each edge of a locally finite non-directed graph is initially
  assigned weight $a>0$. With each step, the particle jumps to a nearest-neighbor vertex.
  The probability of the jump equals to the fraction of the weight attached to the traversed edge in 
  the total sum of the weights of the edges coming out of the vertex where the particle is located
  prior to the jump.
  Each time an edge is traversed, its weight is increased by $1$. In other words, the linear edge-reinforced random walk
  is a random walk on a weighted graph, where the weight of an edge is increased by one each time it is being traversed.

  One of the most important open problems in the theory of reinforced random walks is that of 
  checking if the linear edge-reinforced random walk is recurrent on $\mathbb{Z}^d$ for dimensions 
  $d \geq 2$. Linear edge-reinforced random walk is exchangeable making the model an important 
  example for applying the theorem of de Finetti (see \cite{feller}) and its generalizations (see \cite{df}, \cite{d}).
  The history and some of the most important results in reinforced processes can be found in
 \cite{d}, \cite{davis}, \cite{mr}, \cite{k}, \cite{bk}, \cite{pemantle}, \cite{os} and references therein.

  In this paper we will consider a multi-particle modification of the edge-reinforced random walk model
  similar to some of the reinforced processes studied in \cite{os}.
  We let the walker (or particle) in the edge-reinforced random walk model wait
  an independent exponential time with rate one between the jumps.
  So, the walker jumps from a site to a near by site with rates equal to
  the corresponding ratios. Now we are ready
  to define an $n$-point process $\eta_t=\{\eta_1(t), \dots ,\eta_n(t)\}$, where all $n$ particles
  travel along the edges of a graph $G$, jumping from a site to a neighboring site in $S$,
  the set of all sites.
  Now, let $W_t(e_1),..., W_t(e_k)$ be the weights assigned to all $k$ edges $e_1,...,e_k$
  coming out of a given site $v \in S$ at time $t$. Once again, the initial weights are all assigned
  to be equal to $a>0$, i.e. $W_0(e_1)=W_0(e_2)=\dots=W_0(e_k)=a$. If one of the particles,
  say $\eta_j$, is at site
  $v$ at $t$ when its exponential clock rings, then the particle traverses $e_i$ ($1 \leq i \leq k$)
  with the rate $={W_{t-}(e_i) \over W_{t-}(e_1)+ \dots + W_{t-}(e_k)}$. In which case the corresponding edge 
  weight increases by $1$, i.e. $W_t(e_i)=W_{t-}(e_i)+1$. The recurrence/transience 
  questions arising in this more general model
  are as important as the corresponding questions in the theory of one-article edge-reinforced random walks.

  The edge-reinforced process on $\mathbb{Z}$ with drift $\Delta>0$ can be defined in the following
  way: if a particle is at site $v \in \mathbb{Z}$ at the jump time $t$, then the probability of the particle jumping to
  $v+1$ is $${W_{t-}(v,v+1)+\Delta \over W_{t-}(v-1,v)+W_{t-}(v,v+1)+\Delta}$$ while the probability
  of it jumping to $v-1$ is $${W_{t-}(v-1,v) \over W_{t-}(v-1,v)+W_{t-}(v,v+1)+\Delta}.$$ 
   On trees, the edge-reinforced process with a toward-the-root drift $\Delta>0$ can be stated 
   accordingly. 

  We will concentrate on the most basic case of multi-particle reinforced processes: the two point
  reinforced process $\eta_t=\{\eta_1(t),\eta_2(t)\}$ on $\mathbb{Z}$ with drift $\Delta \geq 0$.
  We will describe an unusual coupling construction associated with the process and as a consequence prove 
  the recurrence of ${(\eta_2(t)-\eta_1(t))}$ whenever $0 \leq \Delta <1$.
  In the case of a two-particle process on $\mathbb{Z}$, one of  the two particles is located to the left of another, except for the times when 
  both particles
  are at the same site. We will denote by $l_t$ the location of the left particle, and by $r_t$
  the location of the right particle at time $t$. When $r_t=l_t$ there is no need to distinguish
  between the ``left" and the ``right" particles. The difference becomes apparent only when
  one of the particles leaves the site. So,  $\eta_t=\{l_t,r_t\}$, and here is the main result of this paper:
  \begin{thm} \label{rec}
  For all $0 \leq \Delta <1$ and $a>0$, $(r_t-l_t)$ is recurrent.
  \end{thm}

  Let us begin by reviewing the Polya's urn model. The urn initially contains $R_0$ marbles of red color,
  land $B_0$ marbles of blue color. We fix a positive integer number $D$. A marble is
  randomly and uniformly drawn from the urn, returned, and $D$ marbles of the same color 
  are added.
  Let $R_n$ and $B_n$ be respectively the number of red and blue marbles in the
  urn after $n$ drawings, and let $\rho_n={R_n \over R_n+B_n}$ be the fraction
  of the red marbles in the urn after $n$ drawings. It is easy to show that $\rho_n$
  is a martingale, and therefore, by martingale convergence theorem, converges to
  a random variable. That random variable $\rho_{\infty}$ is
  in turn shown to be a beta random variable with parameters ${R_0 \over D}$
  and ${B_0 \over D}$, i.e. one with beta density function
  \begin{eqnarray} \label{betaUrn}
  {1 \over \beta({R_0 \over D},{B_0 \over D})} x^{{R_0 \over D}-1}(1-x)^{{B_0 \over D}-1},
  \end{eqnarray}
  where $\beta(a,b) = {\Gamma(a+b) \over \Gamma(a) \Gamma(b)}$.
  One can check that the urn model is {\it exchangeable} (see \cite{feller}), that is if one permutes
  the results of $m$ consecutive drawings, the probability of the outcome does not change.
  By de Finetti's theorem, conditioned on $\rho_{\infty}$, the results of the
  drawings are independent Bernoulli trials, where each time a red marble is selected with probability $\rho_{\infty}$
  and a blue marble is selected with probability $1-\rho_{\infty}$.

  The model trivially extends to the case when $R_0,B_0$ and $D$ are positive {\it real}
  numbers, as well as when there are more than two different types of marbles.
  For instance, consider the case when there are three types of marbles, red, blue and green, in the urn.
  If we start with the amounts $R_0$, $B_0$ and $G_0$ of respectively red, blue and green marbles,
  then the limiting fractions vector will be a Dirichlet distributed random vector, i.e. the cumulative
  density function for the limiting fractions of red and blue marbles will be
  \begin{eqnarray*}
  f(x,y)=
      {\Gamma \left( {R_0+B_0+G_0 \over D} \right) \over \Gamma \left( {R_0 \over D} \right)
      \Gamma \left( {B_0 \over D} \right)\Gamma \left( {G_0 \over D} \right)}
      x^{{R_0 \over D}-1}y^{{B_0 \over D}-1}(1-x-y)^{{G_0 \over D}-1}
    \text{ if }x>0, y>0 \text{ and } x+y <1,
  \end{eqnarray*}
  where the above density is derived by applying (\ref{betaUrn}) twice.

  See \cite{feller} for basic facts on exchangeability, the Polya's urn model and
  a simple version of de Finetti's theorem. A simple
  proof of the convergence to beta distribution can be found in \cite{earadhan}.

  Polya's urns were used to study linear edge-reinforced random walks (see \cite{pemantle}). 
  There, if the walk lives on an acyclic graph, say $\mathbb{Z}$, we can assign a Polya's urn 
  for each site. When the walker is at site $v$, we do a drawing from the urn associated with $v$, where
  the number of the red (respectively blue) marbles in the urn
  is equal to the weight attached to the edge $[v-1,v]$ (respectively $[v,v+1]$) at the time.
  If a red marble is drawn, the walker jumps to $v-1$ and we
  add $D=2$ red marbles to the urn associated with $v$. Similarly, if a blue
  marble is drawn, the walker jumps to $v+1$ and we
  add $D=2$ blue marbles. We do so because
  the graph is acyclic: if the particle ever returns to the vertex, it will be from
  the same direction it took when it left the vertex. For example, if the red marble is drawn,
  the walker will traverse the edge $[v-1,v]$ twice before returning to $v$, thus increasing the weight of
  the edge exactly by $D=2$. The initial conditions $[R_0(v),B_0(v)]$ at an urn associated 
  with site $v$ should be set
  equal to the weights attached to the edges $[v-1,v]$ and $[v,v+1]$
  respectively at the time of the first arrival to $v$. As an example, consider
  the edge-reinforced random walk on $\mathbb{Z}$ that begins at site $0$.
  There the correct initial conditions for a Polya's urn assigned to site $v \in \mathbb{Z}$ should be set equal to
  $$[R_0(v),B_0(v)]=\begin{cases}
      [a, a+1] \text{ if } v<0, \\
      [a, a] \text{ if } v=0, \\
      [a+1, a] \text{ if } v>0.
\end{cases}$$
  What follows is that one can do an infinite number of drawings independently for each of the
  Polya's urns associated with the vertices of an acyclic graph before the walk begins,
  thus completely predetermining
  the trajectory of the walker. Now, the exchangeability property of Polya's urns and de Finetti's theorem mentioned above allows 
  one to restate the edge-reinforced random
  walk as a random walk in random environment (RWRE), where the environment is distributed
  as the limiting beta random variables obtained for Polya's urn processes associated with each 
  vertex of
  the acyclic graph. After that, the large deviation and other techniques are of use in answering
  the corresponding recurrence/transience questions for the RWRE model (see \cite{pemantle}).

  Does the same approach work for the two point
  process $\eta_t=\{l_t,r_t\}$ on $\mathbb{Z}$? The answer is ``no". Consider the case when
  the drift $\Delta=0$. Suppose there is an urn at each vertex of $\mathbb{Z}$. Suppose a vertex $v$ is visited by the right 
  particle $r_t$, and the drawing was done from the urn associated with $v$ and a blue marble 
  was selected, so that the right particle $r_t$ jumps to $v+1$.
  We {\it cannot} add $D=2$ blue marbles into the urn, as it could happen that the left
  particle $l_t$ arrives to the urn from $(-\infty, v-1]$ before the right particle $r_t$ returns to $v$
  from $[v+1, \infty)$. 
  In the latter case, there will be more blue marbles than the weight amount attached to the edge
  $[v,v+1]$ on the right and the rates will not agree. In other words, the representation with Polya's urns and
  similar approaches will not work because the two-point linear edge-reinforced process is nonexchangeable. 
  The non-exchangeability of the process was the main obstacle for studying it as well as for proving 
  Theorem \ref{rec}.

\section{The Polya's urn modified} \label{urn}
  Although the representation with classical Polya's urns fails for the two-point process
  $\eta_t=\{l_t,r_t\}$, there is a way to modify it. Suppose that at each vertex of $\mathbb{Z}$, 
  the associated urn contains not only the red and blue marbles, but also a special
  marble, called {\it magic marble}, such that when the left particle $l_t$ arrives to the site, 
  the magic marble becomes red, while when the right particle $r_t$ arrives to the site, 
  the magic marble becomes blue. Each urn will contain exactly one magic marble in addition 
  to red and blue marbles. When magic marble is selected,
  two marbles of the color assumed by the magic marble will be added into the urn.
  Once again the particles move according to the colors of marbles selected from the urns.
  In other words, if the magic marble is selected when it is red, two more red marbles will be added to
  the urn and the particle will jump left. Similarly, if the magic marble is selected when it is blue, 
  two more blue marbles will be added to the urn and the particle will jump right.

  Let $R_t(v)$ and $B_t(v)$ denote respectively the number of red and the number of blue marbles
  inside an urn associated with site $v \in \mathbb{Z}$, at time $t$.
  The initial number of red and blue marbles,  $R_0(v)$ and $B_0(v)$,
  together with the magic marble must represent the corresponding
  weights  assigned to edges $[v-1,v]$ and $[v,v+1]$ at the time of the first arrival to $v$
  by any of the two particles.
  The left particle is the first to visit the sites to the left of $l_0$, i.e. all $v<l_0$, and the right particle is  
  first to visit the sites to the right of $r_0$.  Hence, for all $a>0$,
  the following must be the initial configuration of red and
  blue marbles assigned to the urns associated with sites in $\mathbb{Z}$:
   \begin{eqnarray} \label{initial}
      [R_0(v),B_0(v)]= \begin{cases}
      [a-1,1+a+\Delta] \text{ if } v<l_0, \\
      [a-1,a+\Delta] \text{ if } v=l_0, \\
      [a,a+\Delta] \text{ if } l_0<v<r_0, \\
      [a,a-1+\Delta] \text{ if } v=r_0, \\
      [a+1,a-1+\Delta] \text{ if } r_0<v
      \end{cases}
    \end{eqnarray}
   {\bf plus} a magic marble in every urn. Here we allow $a-1 <0$ since there is also a magic marble in the urn,
   which is red when the left particle is at the site, and blue when the right particle is at the site.

   We now explain the reason why the magic marble was introduced. First we check that the
   above urn representation produces correct rates up until the first recurrence time
   $\tau_1:=min\{t:\text{ } l_t=r_t\}$. We consider the case when the right particle departs from
   site $v$ at jump time $t <\tau_1$. Suppose that the next arrival to $v$ happens before $\tau_1$,
   then there are three possible scenarios.
   
   {\it Case I: the right particle jumps to the left, and returns to $v$ before the left particle
   arrives.}  So $r_{t-}=v$, $r_{t}=v-1$, $l_t<v$, and we need to add two red marbles into the urn,  
   i.e.  $R_t(v)=R_{t-}(v)+2$. The magic marble stays blue, and as it was the case with one particle
   ERRW model, the rates agree.

   {\it Case II: the right particle jumps to the right, but returns to $v$ before the left particle arrives.} 
   That is $r_{t-}=v$, $r_{t}=v+1$, $l_t<v$ and we need to add two blue marbles into the urn, 
   i.e. $B_t(v)=B_{t-}(v)+2$. The rates agree since the right particle returns to $v$ before the next 
   visit to $v$ by the left particle. Again, the magic marble stays blue, and as it was the case with one particle
   ERRW model, the rates agree.
   
   {\it Case III: the right particle jumps to the right, and the left particle arrives to $v$ before the right
   particle returns from $[v+1, +\infty)$.} This is the case where the chameleon  property of
   the magic marble is used. Once again $r_{t-}=v$, $r_{t}=v+1$, $l_t<v$ and we need to add two 
   blue marbles into the urn,  i.e. $B_t(v)=B_{t-}(v)+2$. 
   Before the departure of the right particle from site $v$ at time $t$,
   $$\text{the weight assigned to }[v-1,v]\text{ was }=R_{t-}(v)$$
   and
   $$\text{the weight assigned to }[v,v+1]\text{ was }=B_{t-}(v)+1$$
   as the magic marble was blue in the presence of the right particle. When the left particle jumps to 
   $v$ from $v-1$ at time $t_1 \in (t,\tau_1)$ before the return of the right particle,  the magic marble
   re-colors into red, and
   $$\text{the weight assigned to }[v-1,v]\text{   }=R_{t_1}(v)+1=R_{t-}(v)+1$$
   and
   $$\text{the weight assigned to }[v,v+1]\text{   }=B_{t_1}(v)=B_t(v)=B_{t-}(v)+2.$$
   One can see that the weights are correct since each edge $[v-1,v]$ and $[v,v+1]$ was traversed exactly
   once. 
   
   That explains why adding magic marble works. The case when the left particle is at site $v$ can be checked by 
   the analogy with the case above. 
   
   Observe that the above coupling of the urn process with $\eta_t$ works
   \textbf{only} up until time $\tau_1=min\{t:\text{ } l_t=r_t\}$, the first time that the particles meet.
   Now, we need to show that the particles meet at least once.
   Therefore before proving Theorem \ref{rec}, we will need to prove the following one-time recurrence result:
  \begin{thm} \label{rec1}
  For all $0 \leq \Delta < 1$ and all $a>0$, $\tau_1 < \infty$.
  \end{thm}
  Since the above urn process is coupled with the above described urn process until the decoupling time $\tau_1$,
  it suffices to prove Theorem \ref{rec1} for the urn process.
  Later it will be shown that the construction will also imply the full recurrence,
  i.e. Theorem \ref{rec}.

\section{Recurrence via coupling.} \label{proof}

  The urn construction construction defined in the preceding section determines $\eta_t=(l_t,r_t)$ for $0 \leq t \leq \tau_1$.
  Here we let $q^l_t(v)$ and $p^l_t(v)$ be respectively the left and the right jump rates for
  the {\it left} particle at site $v$.
  We also denote by  $q^r_t(v)$ and $p^r_t(v)$ respectively the left and the right jump rates for
  the {\it right} particle at site $v$.
  By construction,
  $$q^l_t(v)={R_t(v)+1 \over R_t(v)+B_t(v)+1} \quad \text{ and } \quad
      p^l_t(v)={B_t(v) \over R_t(v)+B_t(v)+1},$$
   and similarly,
     $$q^r_t(v)={R_t(v) \over R_t(v)+B_t(v)+1} \quad \text{ and } \quad
      p^r_t(v)={B_t(v)+1 \over R_t(v)+B_t(v)+1}.$$

   We recall that in the case of edge-reinforced random walks on $\mathbb{Z}$,
   the urn drawings predetermined the outcome of the whole process. There
   the results of all drawings from the urns associated with all the vertices of
   the graph determined uniquely the trajectory of the walker.
   In the one-particle case one determines the limiting fractions of blue marbles
   for all sites in the form of  respective independent beta random variables.
   Then one interprets the walk as a birth and
   death chain with these rates. We want to implement a similar trick for the
   two-point edge-reinforced process.
   We will embed the urn process defined in the preceding section into a three-color Polya's urn process.

\subsection{Magic family} \label{another}
  For each site $v$ and the urn associated with $v$, we define the {\it magic family}
  as all the marbles that were added as the result of selecting
  the magic marble, plus the magic marble itself.  Here is how the magic family is constructed:
  at the beginning the magic family consists of only the magic marble itself.
  When the magic marble is selected from the urn for the first time, and
  two marbles (of one of the two colors) are added to the urn, we include the two into the
  magic family. Each time a marble from the magic family is selected and two new marbles
  are added, we let the two marbles into the magic family regardless of their color.

  All the red marbles that are not in the magic family will be called {\bf pure} red, and all the
  blue marbles that are not in the magic family will be called {\bf pure} blue. Observe that
  for each site $v$, the urn associated with $v$ is a Polya's urn with respect to three types
  of marbles: pure red, pure blue and the marbles in the magic family. Each time when a pure red
  marble is selected, two more pure red marbles are added to the urn. Same is true for pure
  blue marbles.

   For each site $v$, we let $\bar{R}_n(v)$, $\bar{B}_n(v)$ and $\bar{M}_n(v)$
   denote respectively the number of pure red marbles, the number of
   pure blue marbles and the number of marbles in the magic family
   inside the urn associated with site $v$ after $n$ drawings.
   The proportion vector of pure red marbles, marbles in the magic family and
   pure blue marbles
   $$\Big[ \bar{R}_n(v), \bar{M}_n(v), \bar{B}_n(v) \Big]  /(\bar{R}_n(v)+\bar{M}_n(v)+\bar{B}_n(v)) $$
   converges to a Dirichlet random vector with parameters
   ${R_0(v) \over 2}$, ${1 \over 2}$ and ${B_0(v) \over 2}$.

    We observe that after $n$ drawings, $\bar{M}_n(v)-1$ marbles in the magic family are
    of either red or blue color. Therefore $\bar{R}_n(v) \leq R_n(v)$ and $\bar{B}_n(v) \leq B_n(v)$.
    Let $\mathcal{B}(\alpha,\beta)$ denote the beta distribution with parameters $\alpha >0$ and $\beta >0$.
    If for each $v$  we define $p^l_{Polya}(v)$ as the limiting fraction of pure blue marbles,
    then $p^l_{Polya}(v)$ will be a beta random variable with parameters
    ${B_0(v) \over 2}$ and ${R_0(v)+1 \over 2}$.
    Looking back at (\ref{initial}), one can write down the corresponding $\mathcal{B}\left({B_0(v) \over 2},
    {R_0(v)+1 \over 2} \right)$ distribution of $p^l_{Polya}(v)$ for each $v \in \mathbb{Z}$.
    For $a \leq 1-\Delta$,
      $$ p_{Polya}^l(v) \text{ is } \begin{cases}
      \mathcal{B}({a+1+\Delta \over 2},{a \over 2})  \text{ if } v<l_0, \\
      \mathcal{B}({a+\Delta \over 2},{a \over 2}) \text{ if } v=l_0, \\
      \mathcal{B}({a+\Delta \over 2},{a+1 \over 2}) \text{ if } l_0<v<r_0, \\
      0 \text{ if }  r_0 \leq v
      \end{cases}$$
   and for $a>1-\Delta$,
      $$ p_{Polya}^l(v) \text{ is } \begin{cases}
      \mathcal{B}({a+1+\Delta \over 2},{a \over 2})  \text{ if } v<l_0, \\
      \mathcal{B}({a+\Delta \over 2},{a \over 2}) \text{ if } v=l_0, \\
      \mathcal{B}({a+\Delta \over 2},{a+1 \over 2}) \text{ if } l_0<v<r_0, \\
      \mathcal{B}({a-1+\Delta \over 2}, {a+1 \over 2}) \text{ if } v=r_0, \\
      \mathcal{B}({a-1+\Delta \over 2},1+{a \over 2}) \text{ if } r_0<v.
      \end{cases}$$
      Similarly, $q^r_{Polya}(v)$ defined as the limiting fraction of pure red marbles in the urn
      will be a beta random variable with parameters
      ${R_0(v) \over 2}$ and ${B_0(v)+1 \over 2}$.
      So for $a \leq 1$,
      $$ q_{Polya}^r(v) \text{ is } \begin{cases}
      0 \text{ if } v \leq l_0, \\
      \mathcal{B}({a \over 2},{a+1+\Delta \over 2}) \text{ if } l_0<v<r_0, \\
      \mathcal{B}({a \over 2},{a+\Delta \over 2}) \text{ if } v=r_0, \\
      \mathcal{B}({a+1 \over 2},{a+\Delta \over 2}) \text{ if } r_0<v.
      \end{cases}$$
    and for $a>1$,
      $$ q_{Polya}^r(v) \text{ is } \begin{cases}
      \mathcal{B}({a-1 \over 2},1+{a+\Delta \over 2}) \text{ if } v<l_0, \\
      \mathcal{B}({a-1\over 2}, {a+1+\Delta \over 2}) \text{ if }  v = l_0, \\
      \mathcal{B}({a \over 2},{a+1+\Delta \over 2}) \text{ if } l_0<v<r_0, \\
      \mathcal{B}({a \over 2},{a+\Delta \over 2}) \text{ if } v=r_0, \\
      \mathcal{B}({a+1 \over 2}, {a+\Delta \over 2}) \text{ if } r_0<v.
      \end{cases}$$
    Let $q_{Polya}^l := 1-p_{Polya}^l$ and $p_{Polya}^r:=1-q_{Polya}^r$.
    Then the pairs $(q_{Polya}^l(v),p_{Polya}^l(v))_{v \in \mathbb{Z}}$ and
    $(q_{Polya}^r(v),p_{Polya}^r(v))_{v \in \mathbb{Z}}$ can be viewed as two dependent
    random environments. We define $l_t^{Polya}$ as a random walk in
    the random environment $(q_{Polya}^l(v),p_{Polya}^l(v))_{v \in \mathbb{Z}}$
    that starts at $l^{Polya}_0=l_0$ and jumps from  $v$ to  $v+1$ with rate $p_{Polya}^l(v)$ 
    or to $v-1$ with rate $q_{Polya}^l(v)$.
    Similarly, we define $r_t^{Polya}$ as a random walk in
    the random environment $(q_{Polya}^r(v),p_{Polya}^r(v))_{v \in \mathbb{Z}}$
    that starts at $r^{Polya}_0=r_0$ and jumps from $v$ to $v+1$ with rate $p_{Polya}^r(v)$ or to $v-1$ with rate
    $q_{Polya}^r(v)$.

    We observe that conditioned on the two dependent environments, the two random
    walks, $l^{Polya}_t$ and $r^{Polya}_t$, can coexist as two independent birth and death chains.
    In the next subsection we will couple $\{l^{Polya}_t,r^{Polya}_t \}$ with $\{l_t,r_t\}$ so that
    $$l^{Polya}_t \leq l_t \leq r_t \leq r^{Polya}_t ~~~\text{ for } 0 \leq t \leq \tau_1~.$$
    Then showing the recurrence of $(r^{Polya}_t-l^{Polya}_t)$ will prove Theorem \ref{rec1}.
    Here is the heuristic explanation for the coupling construction to follow. If in the original
    urn model $\{l_t,r_t\}$, we substitute the magic marble with a red marble in every urn, then the
    left particle process $l_t$ will be distributed as the random walk in random environment $l^{Polya}_t$.
    If in turn we substitute the magic marble with a blue marble in every urn, then the right particle
    process will have the same distribution as $r^{Polya}_t$. Observe that this heuristics can be generalized to 
    work in the case of more than two particles, e.g. three-particle linear edge-reinforced processes on $\mathbb{Z}$.

\subsection{Coupling with RWRE}
    We notice that the process $\{l_t,r_t\}$ can be predetermined by
    the results of all drawings from the Polya's urns with three types of marbles: pure red, pure blue and magic family.
    We can first do the drawings from the urns associated with all the sites in
    $\mathbb{Z}$, determining $[\bar{R}_n(v), \bar{M}_n(v), \bar{B}_n(v)]_{n=0,1,2,...}$ for all $v$
    and the limiting fractions
    $\{q_{Polya}^r(v), p_{Polya}^l(v)\}_{v \in \mathbb{Z}}$ of pure red and pure blue marbles.
    By de Finetti's theorem, conditioned on $q_{Polya}^r(v), p_{Polya}^l(v)$, $[\bar{R}_n(v), \bar{M}_n(v), \bar{B}_n(v)]_{n=0,1,2,...}$
    are determined as independent trials with probabilities $q_{Polya}^r(v)$ and $p_{Polya}^r(v)$ for pure red and pure blue marbles
    respectively.
    When one of the two particles visits site $v$,  the
    results $[\bar{R}_n(v), \bar{M}_n(v), \bar{B}_n(v)]$ of the $n$-th drawing (if it is the time of $n$-th
    departure from the site) determine the destination site and the coloring of the marbles in the
    magic family.

    We recall that in the Polya's urn model the limiting fraction of marbles of one color
    is a particular beta random variable and   the density function $f_{\{q_{Polya}^r(v), p_{Polya}^l(v)\}}(x,y)$ for
    the pair of limiting fractions $q_{Polya}^r(v)$ and $p_{Polya}^l(v)$ is Dirichlet
    $$f_{\{q_{Polya}^r(v), p_{Polya}^l(v)\}}(x,y)=
      {\Gamma \left( {R_0(v)+B_0(v)+1 \over 2} \right) \over \Gamma \left( {R_0(v) \over 2} \right)
      \Gamma \left( {B_0(v) \over 2} \right)\Gamma \left( {1 \over 2} \right)}
      x^{{R_0(v) \over 2}-1}y^{{B_0(v) \over 2}-1}(1-x-y)^{-{1 \over 2}}$$
      if $x>0, y>0$ and $x+y \leq 1$. 
      
    \vskip 0.3 in 
    Now we construct the coupled process $\{ l^{Polya}_t, l_t, r_t, r^{Polya}_t \}$:
    \begin{itemize}
    \item We condition on the Dirichlet variables $\{q_{Polya}^r(v), p_{Polya}^l(v)\}_{v \in \mathbb{Z}}$. 
    \item We begin with $l^{Polya}_0=l_0 < r_0=r^{Polya}_0$. 
    \item The trajectories of $l_t$ and $r_t$ are determined by drawings from the urns.  
    \item $l^{Polya}_t$ and $r^{Polya}_t$ are independent birth and death chains with corresponding probabilities 
    $(q_{Polya}^l(v),p_{Polya}^l(v))_{v \in \mathbb{Z}}$ and $(q_{Polya}^r(v),p_{Polya}^r(v))_{v \in \mathbb{Z}}$ that move 
    independently of each other and of $l_t$ and $r_t$ except for the times when $l^{Polya}_t=l_t$ or $r^{Polya}_t=r_t$. 
    \item When  the left particles $l$ and $l^{Polya}$ happen to be at the
    same site $v$, they wait for the departure time with rate one. Then the results of the next drawing from the corresponding urn
    are studied. If the selected marble is red or magic marble, both particles jump to $v-1$. 
    However, if the selected marble is blue, the left particle $l$ jumps
    to $v+1$, while $l^{Polya}$ jumps to $v+1$ only if the marble is pure blue, and jumps to $v-1$
    otherwise. The probability that $l^{Polya}$ jumps to $v+1$ is equal to $p_{Polya}^l(v)$, and each such
    drawing is independent of all the other drawings. $l^{Polya}_t$ will still be a random walk in random environment
    $(q_{Polya}^l(v),p_{Polya}^l(v))_{v \in \mathbb{Z}}$, while $l^{Polya}_t \leq l_t$ is preserved.
    \item Similarly, when  the right particles $r$ and $r^{Polya}$ happen to be at the
    same site $v$, they wait for the departure time with rate one. Then
    the results of the next drawing are studied. If the selected marble is a blue or a magic marble,
    both particles jump to $v+1$. If the marble is red, the right particle $r$ jumps
    to $v-1$, while $r^{Polya}$ jumps to $v-1$ only if the marble is pure red, and jumps to $v+1$
    otherwise. The probability that $r^{Polya}$ jumps to $v-1$ is equal to $q_{Polya}^r(v)$, and each
    drawing is independent of the others. $r^{Polya}_t$ will still be a random walk in random environment
    $(q_{Polya}^r(v),p_{Polya}^r(v))_{v \in \mathbb{Z}}$, while $r_t \leq r^{Polya}_t$ is preserved.
   \end{itemize}
   
   In the above coupled process, conditioned on the environments, the independence of $l^{Polya}_t$ and $r^{Polya}_t$
   is preserved and $$l^{Polya}_t \leq l_t \leq r_t \leq r^{Polya}_t~.$$
   Thus showing the recurrence of $(r^{Polya}_t-l^{Polya}_t)$ is enough to prove Theorem \ref{rec1}.

  \vskip 0.3 in
  {\it Recall form the theory of RWREs:} a RWRE on $\mathbb{Z}$ with the right jump
  probability $p(v)$ chosen to be $\mathcal{B}({a+1 \over 2}, {a+\Delta \over 2})$
  distributed at all sites $v$ is transient to the right whenever $0 \leq \Delta<1$, explaining
  the bound on the drift $\Delta$ in Theorem \ref{rec}. In general,
  RWRE on $\mathbb{Z}$ is a.s. transient to the right if and only if $E[\log \Big({p(v) \over 1-p(v)}\Big) ]>0$
  (see \cite{solomon} for the proof, and \cite{sinai}, \cite{kesten} and references therein for more on the subject).
  Not surprisingly the RWRE with the environment $\{p(v)\}_{v \in \mathbb{Z}}$ 
  independently $\mathcal{B}(\alpha_1, \alpha_2)$
  distributed with $\alpha_1>\alpha_2$ is a.s. transient to the right:
  \begin{eqnarray*}
  E[\log \Big({p(v) \over 1-p(v)}\Big) ]
  & = & \frac{1}{\beta(\alpha_1, \alpha_2)}
            \int_0^1 \log \Big({x \over 1-x}\Big) x^{\alpha_1-1} (1-x)^{\alpha_2-1} dx \\
  & = & \frac{4}{\beta(\alpha_1, \alpha_2)}
            \int_{-\infty}^{\infty} s  e^{2\alpha_1 s} (1+e^{2s})^{-(\alpha_1+\alpha_2)} ds \\
  & = & \frac{2^{3-\alpha_1-\alpha_2}}{\beta(\alpha_1, \alpha_2)}
            \int_{0}^{\infty} s  \sinh ((\alpha_1-\alpha_2)s) (\cosh s)^{-(\alpha_1+\alpha_2)} ds \\
  & > & 0,
  \end{eqnarray*}
  where we substitute ${x \over 1-x}=e^{2s}$.
  In the above general case,
  $E \Big[{1-p(v) \over p(v)}\Big]= {\beta(\alpha_1-1,\alpha_2+1) \over \beta(\alpha_1, \alpha_2)}
   ={\alpha_2 \over \alpha_1-1}$ if $\alpha_1 >1$,
  and $E \Big[{1-p(v) \over p(v)}\Big]=+\infty$ if $\alpha_1 \leq 1$.
  
  The expected time of return is a.s. finite only when
  $E \Big[{1-p(v) \over p(v)}\Big]<1$ (see \cite{solomon}). That is only when
  $\alpha_1 > 1+\alpha_2$. Now, if one considers a RWRE on $\mathbb{Z}$ with independent right jump
  probabilities $\{p(v)\}_{v \in \mathbb{Z}}$, each $\mathcal{B}({a+1 \over 2}, {a+\Delta \over 2})$ distributed,
  the recurrence time is infinite since ${a+1 \over 2} < 1+ {a+\Delta \over 2}$ for $\Delta>0$.

  The above implies that the expected time of return to site $v$ for $r^{Polya}$ is infinite
  if $v<r^{Polya}_t$ and is finite if $r^{Polya}_t<v$. Similarly, the expected time of return
  to site $v$ for $l^{Polya}$ is infinite if $l^{Polya}_t<v$ and is finite if $v<l^{Polya}_t$.

\subsection{Proof of Theorem \ref{rec1}.}

  In this subsection we will prove Theorem \ref{rec1}.
  
  \vskip 0.3 in
  \noindent
  {\it Proof of Theorem} \ref{rec1}:
   Recall that the environments $\{(p_{Polya}^l(v),q_{Polya}^l(v)\}_{v \in \mathbb{Z}}$
  and\\  $\{(p_{Polya}^r(v),q_{Polya}^r(v)\}_{v \in \mathbb{Z}}$ of $l^{Polya}_t$ and
  $r^{Polya}_t$ are dependent, but conditioned on the environments,
  the walks $l^{Polya}_t$ and $r^{Polya}_t$ are independently.  We claim that, even though
  the expected time of return of $ (r^{Polya}_t-l^{Polya}_t) $ to zero is
  infinite, $ (r^{Polya}_t-l^{Polya}_t) $ is {\bf recurrent}. The problem can be summarized by
  the following more general lemma.
  \begin{lem}
  Let $p_1(1), p_1(2), ...$ be i.i.d. random variables defined on $(0,1)$ with 
  $$\mu_1 = E\left[ \log\left({1-p_1(i) \over p_1(i)}\right)\right]>0,$$
  e.g. $p_1(i) \sim \mathcal{B}(a_1, b_1)$ for $0<a_1 < b_1$, and let $p_2(1), p_2(2),...$ be i.i.d.
  random variables defined on $(0,1)$ with 
  $$\mu_2 = E\left[ \log\left({1-p_2(i) \over p_2(i)}\right)\right]>0,$$
  e.g. $p_2(i) \sim \mathcal{B}(a_2, b_2)$ for $0<a_2 < b_2$.
  Also let $p_1(0)=p_2(0)=1$.

  If $p_1(0), p_1(1), p_1(2), ...$ are the forward rates for the birth-and-death chain $Z^r_t$, and \\
  $p_2(0), p_2(1), p_2(2), ...$ are the forward rates for the birth-and-death chain $Z^r_t$, then the two dimensional 
  RWRE $X_t=(Z^l_t,Z^r_t)$ returns to zero infinitely often.
  \end{lem}
  
  The following is equivalent statement that we can apply in our case:
   Suppose $Z^r_t$ is a RWRE on $\mathbb{Z}_+$ such that for any $i \geq 0$,
  $Z^r_t$ jumps from site $i$ to $i+1$ with rate $p_1(i)$ and to $i-1$ with rate $1-p_1(i)$,
  and suppose $Z^l_t$ is a RWRE on $\mathbb{Z}_-$ such that for any $i \geq 0$,
  $Z^l_t$ jumps from site $-i$ to $-(i+1)$ with rate $p_2(i)$ and to $-(i-1)$ with rate $1-p_2(i)$.
  If conditioned on the environments $\{ p_1(i)\}_i$ and $\{p_2(j)\}_j$, $Z^r_t$ and $Z^l_t$ are 
  independent birth and death chains, then $Z^r_t-Z^l_t$ is recurrent.
  
  The proof of the above lemma can be thought of as an exercise on use of harmonic functions in stochastic
  processes. It can be done with Lyapunov functions (see \cite{fmm}), or alternatively with conductivities. 
  $~~~\square$

\subsection{Proof of Theorem \ref{rec}} \label{a=1}

 It was essential for proving Theorem \ref{rec1} that $\eta_t= (l_t, r_t)$ was defined via urns and magic marbles for $t \leq \tau_1$.
 There are many ways to complete the proof of Theorem \ref{rec},
 one is to notice that when the particles separate after the first meeting time
 $\tau_1$, we can do the whole coupling construction (that lead us to the proof
 of Theorem \ref{rec1}) anew, starting from scratch. The only thing different will be the
 initial marble configuration for $\eta_t=(l_t, r_t)$,
 and the environments of $l^{Polya}_t$ and $r^{Polya}_t$, but only at finitely
 many sites, thus establishing the finiteness of the second meeting time $\tau_2$.
 The finiteness of $\tau_3, \tau_4, \dots$ follows by induction. $~~~\square$
 
 This proof followed from the unusual coupling construction with the magic marbles and the
 domination by two RWREs, one on the right and one on the left. 
 As we already mentioned, the above domination can be constructed without using the magic marble approach,
 and generalized to work for more particles than two. However the approach taken in this paper allows us better understand
 the dynamics behind the two-point processes, and is valuable as an innovative coupling technique.

\section*{Acknowledgment}
 The author wishes to thank all members of UCLA probability group for
 their support in the three years that the author worked there as a postdoc, and Silke Rolles who worked at UCLA at that time 
 for useful discussions. While working on this version of the paper, the author received a lot of encouragement from members of
 the probability group at Oregon State University. The author wishes to thank Robert Burton and Mina Ossiander for
 sharing thoughts on nonexchangeable processes.

\bibliographystyle{amsplain}

\end{document}